# LATTICE BASIS AND ENTROPY

*Vinod Kumar.P.B[1] and K.Babu Joseph[2]*
*Dept. of Mathematics*
*Rajagiri School of Engineering & Technology*
*Rajagiri Valley.P.O, Cochin – 682039*
*Kerala, India.*

**ABSTRACT:**

We introduce the concept of basis for a lattice. This basis plays a vital role to determine the completeness and consistency of the lattice. Weighted lattices are introduced and its complexity is formulated. Some axiomatic systems, considered as lattices, are also studied.

2000 AMS classification: 06B20; 06B23; 06E99.

Key words & Phrases: Lattice, Basis, Entropy etc.

## 1. Introduction

We consider lattices as algebraic system $(L, *, +)$ and as poset $(L, \leq)$ according to the context. We follow the standard definitions as in [1]. In the second section, we define the basis for a lattice and some results are obtained. Logic, Set theory, Geometry and Number theory are lattices as they are axiomatic systems. For all the axiomatic systems the set of axioms is a basis. We discuss the consistency and completeness of these systems from lattice theoretic perception in the third section, The topological entropy as defined by Adler [3] can be interpreted as a measure of complexity of a function in a topological space. For 'nice' interval maps there is a graph theoretic approach to calculate the entropy. Using this approach we introduce complexity in weighted lattices in the fourth section.

1 vinod_kumar@rajagiritech.ac.in
2 babu_65@hotmail.com

## 2. Lattice dimension

***2.1: Definition***: *Let $(L, \leq)$ be a lattice. Let $B \subset L - \{0\}$, where 0 is the least element. $B$ is called a basis for the lattice L if,*

*(i) $x, y \in B \Rightarrow x \leq y, y \leq x$.*

*(ii) If $z \in L$ then there exist $x \in B$ such that $x \leq z$.*

$|B|$ is called dimension (denoted by dim(L)) of the lattice. Elements of $B$ may be called 'basons'.

***2.2: Examples***:

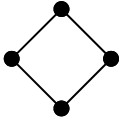

Fig.1

1. Dimension of the above lattice is 2.
2. Dimension of $(N, \leq)$ is 1.
3. $(R, \leq)$ do not have a basis.

***2.3: Proposition***: *Basis of a lattice, if it exists, is unique*.

Proof: If possible, let $B_1$ and $B_2$ be two bases of the lattice L.

Let $x \in B_1 \Rightarrow$ there exists $y \in B_2$ such that $y \leq x$ (since $B_2$ is a basis).

Since $B_1$ is a basis there exists $z \in B_1$ such that $z \leq y$.

$\therefore z \leq y \leq x \Rightarrow z \leq x$, which is a contradiction to the fact that $B_1$ is a basis.

So the basis, if it exists, is unique. •

***2.4: Proposition***: *If L is a chain and if dimension of L exists, then dim (L) =1*.

Proof: Let B be the basis for L with $|B| = n \, (\geq 2)$.

i.e. there are at least n non comparable elements in L. which contradicts the fact that L is a chain. So, $|B| = 1$. $\therefore$ dim (L) =1. •

***2.5: Remarks***: *L is a chain and $0 \in L$ does not imply that there is a basis for L*.

Let $L = [0,1]$ with usual order.

***2.6:*** *Converse of above 2.4 is not true*. dim(L) =1 does not imply that L is a chain.

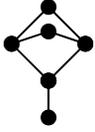

Fig.2.

For example, dimension of the above lattice, which is not a chain, is 1.

***2.7:Proposition***: *If f: $L_1 \to L_2$ is an on to isomorphism and B is a basis for $L_1$ then f(B) is a basis for $L_2$*.

Proof. Let x,y ∈ f(B) $\Rightarrow f^{-1}(x), f^{-1}(y) \in B$

$\Rightarrow f^{-1}(x)$ and $f^{-1}(y)$ are not comparable in $L_1$( since B is the basis for $L_1$).

$\Rightarrow$ x and y are not comparable in $L_2$ ( since f is an isomorphism).

b∈$L_2$ $\Rightarrow$ b = f (a) , for some a ∈$L_1$.

$\Rightarrow a = f^{-1}(b)$

$\Rightarrow$ there exists some x ∈ B such that $x \leq f^{-1}(b)$

$\Rightarrow f(x) \leq b$

Let y = f(x).

So there exists some y ∈ f (B) such that y ≤ b.

∴ f (B) is a basis for $L_2$.•

***2.8: Corollary***: *Two isomorphic lattices have the same dimension*.

***2.9:*** *Converse of 2.8 is not true*.

$L_1$        $L_2$

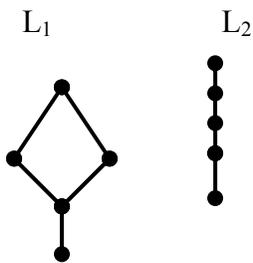

Fig.3.        Fig.4.

dim($L_1$) = dim($L_2$) = 1. But they are not isomorphic.

***2.10: Proposition***: *All basons are atoms. But the converse is not true.*

Proof: If x is a bason then all other elements of the lattice, if they are comparable to x, will be greater than or equal to x. So x can not be written as union join of two elements .•

To show that converse is not true let L be the following lattice.

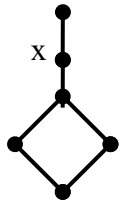

Fig 5

x is an atom, not a bason.

***2.11***: ***Proposition***: *If L is a finite dimensional Boolean algebra then all the atoms of L are also basons*.

Proof: Let the dimension of L be n. L is isomorphic to $(P(X), \subseteq)$ where $|X| = n$. Let f be an isomorphism from L onto P(X). Let x be an atom Then f(x) is a singleton set (If f(x) is not a singleton set then f(x) is union of singletons, hence f(x) can't be an atom).So f(x) is a bason in P(X). Hence by proposition 2.7, x is a bason. .•

***2.12***: ***Proposition***:  A lattice is a finite dimensional Boolean Algebra iff any $x \in L-\{0\}$ can be written as join of basis elements.

Proof: It is obvious from the above proposition and proposition 2.7.•

***2.13:*** Both the implications of Proposition 2.12 is not true in general lattices.

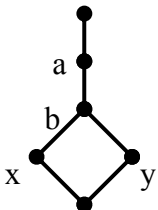

Fig.6

In the above lattice , the basis is B = {x,y}. a ≠ x + y.

*2.14: Definition*: A basis of a lattice is called an orthogonal basis if ∀x ∈ B, x′ ∈ B.

*2.15: Examples*:

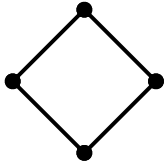 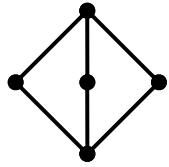 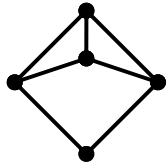 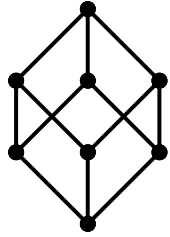

Fig.7.          Fig.8.          Fig.9.          Fig.10.

In the above four lattices the first two (Fig.7 & 8) have orthogonal basis. The basis of last two (Fig.9 &10) are not orthogonal.

*2.16: Proposition:* *If a lattice L with unique complements has an orthogonal basis then dimension of L is even.*

Proof: For every bason x , x′ is also a bason and there are no other complements for x. So, number of basons is even. Hence dimension of L is even.•

*2.17: Definition:* Let L be a bounded lattice. x ∈ L is called an isolated element if x ∗ y =0 and x + y =1 , ∀y ∈L -{0}.i.e. x is isolated if every element(other than 0) are complements of x.

*2.18: Examples:*

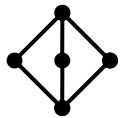 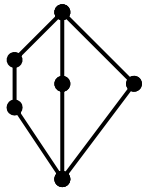

Fig.11.          Fig.12.

In Fig.11 all the elements (except 0 and 1) are isolated. In Fig.12 only the middle one is isolated.

*2.19: Definition:* A lattice without isolated elements is called a consistent lattice, otherwise it is inconsistent.

*2.20: Definition:* if all the basons of a lattice are isolated, then such a basis is called an isolated basis.

*2.21: Proposition:* *Two dimensional Boolean algebra is the only finite dimensional Boolean algebra which is inconsistent.*

Proof: In two dimensional Boolean algebra singletons are inconsistent. If dimension is greater than two there are no inconsistent elements. •

*2.22: Proposition*: *Two dimensional Boolean algebra is the only finite dimensional Boolean algebra which has an orthogonal basis.*

Proof: If there are only two elements in X, complement of singletons are singletons. So the basis is orthogonal. If there are more than two elements in X, complements of singletons are not singletons. •

*2.23: Definition:* Let L be a lattice and $L_1 \subseteq L$. $L_1$ is said to be an independent set if there exist at least one element in $L_1$ which can be written as join of other elements of $L_1$. A set is independent if it is not dependent.

*2.24: Proposition:* *Basis of a lattice is an independent set.*

Proof: Obvious from the definitions. •

*2.25: Proposition:* *A set containing an isolated element is an independent set.*

Proof : If there is an isolated element x in $L_1$, then none of the elements of $L_1$ can be written as join of other elements of $L_1$, because x + y =1, $\forall$y. So $L_1$ is an independent set. •

### 3. Some Axiomatic Systems

If we consider Logic as a lattice of statements bounded by tautology and contradiction then it is inconsistent because there are statements which can be neither derived from any of the other statements and nor any other statements can be derived from it. We usually call them 'conjectures'. Or in other words, conjectures are isolated elements of the lattice of statements. This type of inconsistency is equivalent to Gődel incompleteness.i.e if we prove that there exist atleast one isolated element in a lattice then it is equivalent to proving a generalization of Gődel incompleteness theorem.

*inconsistency* $\Leftrightarrow$ *Godel incompleteness* .

In any axiomatic system, axioms are basons. If there is atleast one isolated element in that system then it is Gődel incomplete. As we increase number of propositions then the complexity involved to reach a conclusion also increases. If we wish to prove a theorem (all theorems belong to the lattice of statements) then in fact we have to start from a set of axioms. From these axioms we prove some propositions. We use this propositions to prove the required theorem. So there is some complexity to prove a theorem. Ultimately this complexity can be regarded as a measure of the lattice. The lattice dimension may be finite, but the lattice can be infinite. For example Euclidean geometry, Set theory, Peano's Arithmetic etc are examples of finite dimensional lattices.

Euclidean geometry is a finite dimensional lattice. The whole geometry is based on 5 postulates of Euclid [2]. All the propositions can be derived from these 5 axioms (of course with some defined concepts). There are a lot of conjectures in Euclidean geometry. The fifth postulate is the most controversial one. Fifth postulate and its negation are consistent with all other axioms. But Geometry is Gődel incomplete.

Set theory is also finite dimensional [4]. There are many conjectures in Set theory. Axiom of choice can be regarded as the most controversial axiom of set theory.

So in a sense, it can be equated with the fifth postulate of Geometry.

Arithmetic is also five dimensional since it is derived from five Peano axioms. Since there are isolated elements it is also Gődel incomplete. If we negate successor axiom, then a new arithmetic will emerge. In short,

*Fifth Postulate* $\Leftrightarrow$ *Axiom of choice* $\Leftrightarrow$ *Successor axiom*.

So naturally a question arises – In every inconsistent lattice does there exist a bason which is controversial (i.e complement of that bason together with all other basons will form a basis for some other lattice)? It is also a factor which determines the complexity of the lattice. Not only that the inconsistency is because of such basons, but also isolated elements are basons. i.e conjectures can be taken as axioms, but the Gődel incompleteness will remain because there are other conjectures. So if a lattice is inconsistent (Gődel incomplete) it will remain as inconsistent even if we add new basons. A paradoxical conjecture is 'how many conjectures will remain as a conjecture in an axiomatic system?'. This question itself is a conjecture in any system. Because inconsistent system will remain as it is it is certain that there will be at least one conjecture. Since all the conjectures are basons and if there are infinitely many conjectures (which remain as conjectures for ever) then the system is infinite dimensional. In this direction it can be conjectured that (1) Set theory is an infinite dimensional axiomatic system (2) Geometry is an infinite dimensional axiomatic system (2) Arithmetic is an infinite dimensional axiomatic system.

**4. Lattice complexity**

If the arguments are vague (fuzzy) then naturally conclusion is uncertain. In the lattice of statements, if we admit this type of fuzzy arguments then such lattices are called weighted lattices.

*4.1: Definition*: A Lattice $(L, \leq)$ is called a weighted lattice if

(i) $x \leq x, \forall x \in L$

(ii) There exist positive real numbers p and q, $p + q \leq 1$ such that

$x \leq_p y \Rightarrow y \leq_q x$.

(iii) $x \leq_p y, y \leq_q z \Rightarrow x \leq_r z$, where $r = \text{Min}(p,q)$.

*4.2: Example:* Consider the Lattice $(L, \leq)$. $L = \{0, x, y, 1\}$ and

$0 \leq_{0.2} x$; $0 \leq_{0.3} y$; $x \leq_{0.5} 1$; $y \leq_{0.5} 1$

$1 \leq_{0.2} x$; $1 \leq_{0.1} y$; $x \leq_{0.5} 0$; $y \leq_{0.7} 0$.

This defines a directed graph with vertex set $\{0, x, y, 1\}$ and Edge set $\{(0,x),(0,y),(x,1),(y,1)\}$ with edge weight as defined above.

From this directed graph we can calculate the entropy as defined in [4] as a measure of complexity.

All the usual lattices are weighted lattices with either p=0, q=1 or q=0, p=1

*4.3: Definition:* Entropy of a weighted lattice $(L, \leq)$ is $\log(\text{Max}(|\lambda|))$ where $\lambda$ is the spectral radius of the adjacency matrix of the directed graph of L. It is denoted by Ent(L).

*4.4: Definition:* A lattice is chaotic if Ent (L) > 0.

*4.5: Proposition: Usual lattices are not chaotic.*

Proof: Since the entries of the adjaceny matrix are either 0 or 1, the matrix will be upper triangular (or lower triangular). So the entropy will be 0. •

So in a sense there is complexity iff arguments are vague.

**Acknowledgement**

The authors thank the management of Rajagiri School of Engineering and Technology, India for the immense support and inspiration given to them.

<p align="center">*********************</p>